# The Birth of Calculus: Towards a More Leibnizian View

Nicholas Kollerstrom

nkastro3@gmail.com

*We re-evaluate the great Leibniz-Newton calculus debate, exactly three hundred years after it culminated, in 1712. We reflect upon the concept of invention, and to what extent there were indeed two independent inventors of this new mathematical method. We are to a considerable extent agreeing with the mathematics historians Tom Whiteside in the 20th century and Augustus de Morgan in the 19th. By way of introduction we recall two apposite quotations:*

*"After two and a half centuries the Newton-Leibniz disputes continue to inflame the passions. Only the very learned (or the very foolish) dare to enter this great killing-ground of the history of ideas" from Stephen Shapin[1] and*

*"When de l'Hôpital, in 1696, published at Paris a treatise so systematic, and so much resembling one of modern times, that it might be used even now, he could find nothing English to quote, except a slight treatise of Craig on quadratures, published in 1693" from Augustus de Morgan[2].*

## Introduction

The birth of calculus was experienced as a gradual transition from geometrical to algebraic modes of reasoning, sealing the victory of algebra over geometry around the dawn of the 18th century. 'Quadrature' or the integral calculus had developed first: Kepler had computed how much wine was laid down in his wine-cellar by determining the volume of a wine-barrel, in 1615,[1] which marks a kind of beginning for that calculus. The newly-developing realm of infinitesimal problems was pursued simultaneously in France, Italy and England. Concerning what 'it' was that was being discovered, Tom Whiteside found that it was tempting 'to admit two criteria into a working definition (without excluding others); first, that differentiation and integration be seen as inverse procedures; and secondly that both be defined with respect to an adequate algorithmic technique.'[2] In addition to these we require what I suggest is Leibniz' criterion, that it be a process whereby 'the imagination is freed from a perpetual reference to diagrams;'[3] i.e., the method is non-visual.

The epochal invention of the differential calculus happened in the latter half of the 17th century. A degree of consensus exists that its appearance has to be located somewhere between 1669, when Newton's manuscript *De Analysi* first appears, and 1696 when de L'Hôpital's textbook on the differential calculus appeared - his *L'Analyse d'Infiniment Petits*. Any priority claim has to be based on a publication or witnessed manuscript *before* that final date. The vexed question is asked, did one or two people - or perhaps more - create it? We are now at the 300th anniversary of the publication of the highly biased Royal Society's own report on the subject,[4] the culmination of what can safely be called the greatest row ever in the history of science.

---

[1] Johannes Kepler, *Stereometrica Doliorum,* Prague 1615 (Solid Geometry of Wine Barrels).
[2] D.T. Whiteside 1961 'Patterns of Mathematical thought in the later 17th Century', *Archive for History of exact Sciences* (AHES) 1:179-388: p.366.
[3] J.M. Childs, *The Early Mathematical Manuscripts of Leibniz*, 1920: translation of *Historia et Origo Calculi Differentialis* composed by Leibniz 1713-5, published 1846, p.25.

Secretly written by Isaac Newton himself, that report dismissed Leibniz as a plagiarist and averred that 'Second inventors have no right,'[5] putting Newton as the first inventor on the grounds that his new method was contained in early manuscripts of his dating back to the 1670s. At the time of making that accusation Newton had indeed published his method, first in *De Quadratura* of 1704 and then an earlier rudimentary version in *De Analysi* of 1711 – but this was all rather late in the day. A Newtonian method using 'moments' with dots on the top of letters had been published in 1693, but not under his name: John Wallis wrote it, as we'll see. Meanwhile Leibniz had published his method in 1684 entitled *A New Method for maxima and minima .... and a Remarkable type of Calculus for them*,[6] that being the seed from which the mighty Continental use of the new method was blossoming. Some take the view that Newton had indeed published on the topic, in his *Principia* of 1687, after all did it not have what retrospectively came to be called a 'fluxions lemma'? We shall in due course have occasion not to accept that view.

## Genesis

By way of attempting to avoid the great errors committed by earlier scholars on this topic, we here cite four more modern examples, of invention's first-beginning, from the history of technology. Confusion can only be avoided if the concept of invention is first properly defined. An inventor has to have performed a *public* act, which *alters* the world by becoming the *seed* from which future growth takes place. The innovation can often be recognised by the sense of delight which the historical characters display around that new-beginning. Briefly therefore we ask the reader's indulgence to do little more than mention these four birth-moments: the steam engine 1804-1814, the motor-car 1862, the electric light 1879, and the laser beam 1960.

1. Invention is an *inherently public act* that involves *risk* and *commitment*. Thus, Richard Trevithic put a steam-engine onto a railway line in Cornwall in 1804. Had he merely kept it in his back-garden he would not be honoured today as an inventor. But, several persons took credit for that invention and George Stevenson probably built the first decent, workable steam engines starting in 1814. Usually but not always notable inventions are associated with one distinctive moment.

2. There are today cars on the street. If one could move backwards in time, there would be fewer and fewer cars – until, reaching 1862, there would be but one car on one street, and that in Paris and that driven by a M. Etienne Lenoir. That moment has to signify the *invention* of the motor-car and it is unique, there was no other. Therefore, M. Lenoir is awarded the credit.

3. Edison made his electric light work on 23 October 1879. Upon hearing that Joseph Swann in Newcastle was claiming to be the inventor of the filament electric light, he commenced a lawsuit: this we might call a 'Newtonian' attitude, whereby 'Second inventors have no right' as Isaac Newton pronounced severely in 1715, concerning the calculus dispute. But, the year

---

[4] Newton 'Account of the Commercium Epistolicum,' *Philosophical Transactions* (PT) January 1714/5, 29:173-224: reproduced in Hall, R. *Philosophers at War,* CUP 1980, 263-314. It appeared towards the end of 1712.
[5] *Account of the Comm. Epist*. P.215: Hall, 1980, p.305.
[6] Leibniz, *Acta Eruditorum,* 1846, 467-73.

1880 saw *both* Britain's House of Commons lit up by Swann's bulbs *and* Edison's Menlo Park laboratory lit up at night. Swann had publicly demonstrated his light-bulb on 18th January 1879. A happy resolution was reached with the setting up on the Edison-Swann light-bulb company, a case, alas rather rare, where synchronous and independent invention is amicably acknowledged.

4. Coherent light was first envisaged on 26 April, 1951 by Charles Townes, while sitting on a park bench in Washington DC. That was a private experience, a eureka-moment, which was embodied, i.e. made to work, a few years later when a laser beam was produced by Theodore Maiman. That happened in his Malibu laboratory, California, on 15th May, 1960 when a beam of coherent light shone forth from a ruby crystal. The new invention causes excitement, which helps the historian in locating the decisive moment. It has to be *witnessed, dateable* (preferably, though this is not essential) and *public*.

## Finding Pi

Leibniz while being taught maths by Christiaan Huygens found his beautiful but slowly-converging circle-squaring formula in Paris in late 1673: $1 - 1/3 + 1/5 - 1/7 + \ldots$ which as he put it gave 'the area of a circle whose diameter is unity.' We express that infinite sum as $\pi/4$. He had no idea that Isaac Newton had earlier invented the obscure circle squaring infinite-series $2RB - B^3/3R - B^5/20R^3 - B^7/56R^5 - 5B^9/576\,R^7$ 'etc'. No-one received any explanation for this cryptic formula until some years later, when in October 1676 Newton wrote to Leibniz. John Collins (who had some access to Newton's manuscripts) had sent the formula to David Gregory on 24 March 1670, with no explanation (*The Correspondence of Isaac Newton*, Vol.I, p.29) then David Gregory reprinted it in a 1684 tract (*Exercitatio Geometrica*

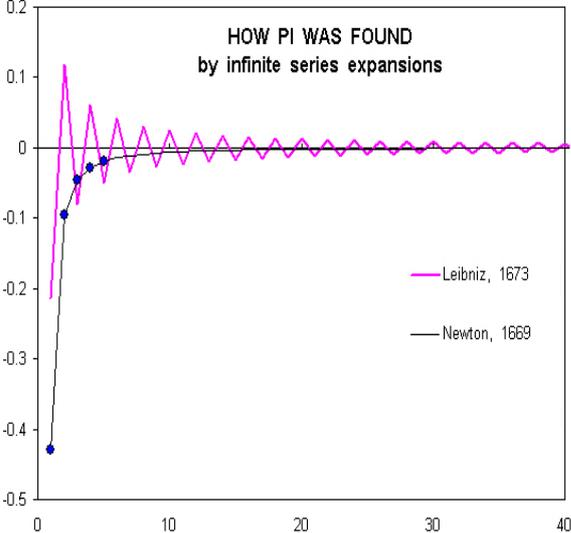

Edinburgh 1684, p.3). Years later, Leibniz did concede Newton's priority in this regard.[7]

To fathom that cryptic formula which Newton came up with, *we today* would need to use the Leibniz-Bernoulli terminology, something like $\int \sqrt{(R^2 - y^2)}\,dx = 2RB - B^3/3R$ etc. where putting $B = R$ gives area of semicircle: Newton had *integrated* (as we would say) the binomial expansion which he had found, of that square-root term. These five terms gave a $\pi$ value of 3.18, so one would require somewhat more to be useful (see figure). Newton's innovation here, was that his act of 'quadrature' had generated an infinite series. What was then called *quadrature* viz.

squaring the circle did then later in the century morph into meaning *integration,* the term coined by Jacob Bernoulli in 1690. By the time Newton published his *De Quadratura* of 1704 it had lost that earlier meaning and that opus despite its title did *not* give any circle-squaring formula.

In 1671 James Gregory wrote out the infinite-series expansion:
$$\text{Arctan } x = x - x^3/3 + x^5/5 - \ldots$$

---

[7] Childs op. Cit. 1846, p.46.

And putting x = 1 would have given Leibniz' series: it has been alluded to as the 'Gregory-Leibniz series' for circle quadrature ever since. Even Leibniz conceded, 'as afterwards it became known it had been worked out by Mr Gregory.'[8] We are by no means in accord with Leibniz here, in conceding the priority of his youthful discovery to another! Maybe Gregory ought to have seen this, but there is no record that he did so. A tangent series is not a circle quadrature. It may be easy in retrospect to say that x = 1 here gives a $\pi/4$ value of radian measure (i.e., the angle whose tangent is one), but they then did not use or speak of $\pi$ in such a manner, and it may not be our business to conjecture as to what historical characters should have done: *that is Leibniz' series*, and not anyone else's. Gregory was subsequently appointed to the chair of mathematics at Edinburgh, then became blind and died in 1675, so maybe had little leisure to develop this theme. Years later, in 1699, Abraham Sharp became the first to use Gregory's tangent formula for obtaining the value of $\pi$, employing it to derive the series for 30° so that it converged more quickly, summing to $\pi/6$ instead of $\pi/4$, suggesting that such an application may have been less evident than was later assumed.

While appreciating Leibniz' modesty we are compelled by a desire for historical accuracy and our definition of 'invention' to reject his judgement over priority. 1676 was the first year when anyone heard of a circle-squaring procedure or formula from Newton, which was three years after Christiaan Huygens and Leibniz understood the circle-squaring convergence. In the mid-1670s the popular play 'The Virtuoso' by Thomas Shadwell satirised the Royal Society, scoffing at its attempts to 'square the circle' as futile and impossible. And yet that immortal dream of philosophers had been realised, in Paris in 1673.[9]

## *De Analysi*

The '*De Analysi*' composed by Newton in 1699 resulted in Newton acquiring the Lucasian chair of mathematics in that year, as a consequence of him presenting it to his teacher Isaac Barrow. It remained unpublished, as he refused Barrow's entreaties to publish it, and likewise responded to John Collin's requests to distribute it with morose comments that he would have it returned and regretted having shown it to anyone. The text did indeed present his infinite binomial expansion method, but without any proof - and hinted that only in the mind of the Deity, not subject to human limitations, could its sum really be known: 'we, mere men possessed of only finite intelligence, can neither designate all their terms [of his infinite equations] nor so grasp them as to ascertain exactly the quantities we desire from them.'[10] Also, it gave a method for finding tangents comparable to that described soon after by John Wallis in the *Phil Trans* of 1672.[11] (adjacent to the article reporting Newton's new reflecting telescope: its tangent diagrams became immortalized by bearing in their midst an image of the new telescope) Gaining the Lucasian professorship, Newton may not have wanted the reputation of one who divided both sides of an equation by zero, for that is what these methods did. His intuition may have told him that this was a workable method – but, he replied with a firm negative to Collins' second request to publish it: 'I wish I could retract what has been done…' He further added mysteriously, concerning *De Analysi,* 'it's plain to me by ye fountain I draw it from, though I will not undertake to prove it to others.' (Newton to Collins, Nov 8, 1676) These things sound rather personal, and unripe for publication.

---

[8] Leibniz, 'Historia et Origo Calculi Differentialis' in Childs 1920, 22-57,46. As Leibniz' riposte to the *Commercium Epistolicum,* (p.5), it was composed 1714-16.
[9] According to Miklos Horvath (*On the Leibnizian Quadrature of the Circle,* Budapest, 1982 (online) pp.75-83,76.) Leibniz found his quadrature in the summer of 1673 and told Huygens about it in December, then alluded to it in correspondence in the summer of 1674, sending an account of it to Huygens in October 1674; publishing it in 1682.
[10] *De Analysi* 1669, MP II 206-247, 243.
[11] *Phil. Trans*. 1672, 7,*Binae Methodi Tangentium* Doctoris Johannis Wallisi, 4010-4016.

When four decades later, that text became cited in the priority dispute, the protagonists seemed unconcerned that Newton had not agreed to publish or distribute it, or that Wallis in both editions of his *Algebra* had discerned therein no hint of any new method of calculus. Much has been made to hinge upon this document, so let us quote Louis Trenchard More's carefully-expressed view, in his biography of Newton:

> *There is much difference of opinion as to whether this work [De Analysi] contained anything on fluxions or the calculus. Newton's partisan's hold that it did, and that an able mathematician could have worked out his secret from it. Brewster, who is sufficiently biased, is doubtful, but such critics as De Morgan, Professor Child and the continental mathematicians, hold that it did not... The methods employed are merely a correlation of what he has learned from Descartes, Wallis and Barrow, combined with his original methods of infinite series and their reversion. In the true sense of the word, he could not differentiate, except as we now say, from first principles.[12]*

For comparison, in 1671 Leibniz submitted his *Theoria Motus Abstracti* to the French Academy (*Op Leibn*. ii part ii p.35), throughout which there occurs a frequent approximation to the idea of infinitely small quantities having a ratio to each other. De Morgan commented on this: 'In 1671 it was working in Leibniz's mind that in the doctrine of infinitely small quantities lay the true foundation of that approach to the differential calculus which Cavalieri presented'.[13]

## Sluse's Method

In January 1673 a Mr Renatus Sluse (Canon of Liege) published a mathematical device whereby a tangent could be fitted to any polynomial curve,[14] Sluse's method derived from that of Hudde, characterised by Whiteside as 'complex' and 'cumbrous.'[15] But, a mere month before its publication, Newton strangely wrote to John Collins, a mathematical fellow of the Royal Society with this method as if it were his – and, as De Morgan observed, it unequivocally was Sluse's method there described[16] – for finding a tangent, as if he were trying to claim priority. Some years later in 1676 Leibniz was shown this letter, and it was averred that he might have learnt various things from it. Then in October 1676 Newton wrote to Leibniz also concerning the Sluse method, making the odd comment: 'Nobody, if he possessed my basis, could draw tangents in any other way, unless he were deliberately

---

wandering from the straight path.'[17] Of this letter, Newton in 1712 affirmed (commenting anonymously in the Royal Society's 1712 *Commercium Epistolicum* judgement) '… which puts it past all dispute that he [Newton] had invented the Method of Fluxions before that time [i.e., October 24th 1676].'[18] But, a person who had the method of calculus, as we would now call it, would neither want nor need this Sluse method, which had no explanation but merely worked: it was a *precursor* of the calculus. Leibniz found it to be a stimulus in discovering his calculus method as his biographer Hoffman showed.[19] Newton later argued that it was by reading his two letters about Sluse's method that Leibniz had gained access to his differential theory.

The Editor of Volume I of the *Correspondence* would seem to agree with the Royal Society's verdict of 1713 concerning that 1672 letter, for he wrote of Newton's 'Sluse' method 'the method is that of a first-order partial differential of a polynomial $\Phi(x,y)$'. (Corr. Vol 1, p.254) Had such been present, it would indeed have settled the matter. Westfall, in his *Never at Rest* repeatedly affirms that the seeds of mistrust and suspicion were sown initially by Leibniz failing to acknowledge Newton's work in his initial calculus publication: 'As we have seen, the storm that burst in 1711 had been gathering for many years, indeed from the moment in 1684 when Leibniz chose to publish his calculus without mentioning what he knew of Newton's progress along similar lines;' (p.742) or, 'nothing could vindicate Leibniz from the charge on which so much of the dispute hinged, that he failed to mention the correspondence of 1676 when he first published.' (p776). The reader is invited to evaluate whether there was any such obligation, and if so what could have been acknowledged.[20] Leibniz did not reach his breakthrough by pondering the series expansions Newton had shown him (as Hoffman has shown, *Leibniz in Paris* 1974), nor did he 'square the circle' via the equations of Gregory, as the English felt he should have done, but by a different method. There might have been less controversy had he made these acknowledgements, but that he was obliged to do so is far from self-evident.

The modern reader, perusing these letters of 1672 and 1676, may tend to agree with Leibniz:

> *'Il n'y a pas la moindre trace ni ombre du Calcul de Differences ou Fluxions dans toutes les anciennes Lettres de M. Newton que j'ai vues…'*[21]

These letters concerned the great discovery Newton had made, that the binomial theorem expands into an infinite series if fractional indices are involved, and how this could be used to 'square the circle'.[22] The mathematics-lecturer and historian Augustus De Morgan wrote of these two letters: 'We know that Newton's letters did not treat of fluxions, nor contain anything from which the writer of a system could draw his materials'[23]. Concerning the 'second letter' of Newton, Hall commented: 'The second letter contained many mathematical treasures, but not the concept of fluxion, nor one example of an expression involving fluxions.[24] It gave special examples of areas found under curves, while hinting that he was in possession of a more general method. The Royal Society's verdict, whereby such vague hints allegedly communicated information to Leibniz, well merits the Shakespearean comment of de Morgan: 'When Glendower said, "I can call spirits from the vasty deep,"[25] no one ever supposed that he "partly described" the "method" of doing it' (De Morgan *Essays* p.157).

---

[17] *The Correspondence of Isaac Newton*, Ed H.Turnbull et. al., CUP, 1959-77, II, Newton to Oldenburg, 24 October 1676: 130-149, 134.
[18] 'An Account of *Comm. Epist.*, p.193; Hall, 1980, p.283.
[19] Hoffman 1974, pp.73-4.
[20] In July 1684 Leibniz wrote to his friend Mencke, about his new calculus article, saying: 'I do not believe either, that Mr Newton will claim it for himself, but only some inventions relating to infinite series which he has in part also applied to the circle.' (Jason Bardi, *The Calculus Wars,* 2006, p.117.)
[21] *Corr.* VI, 1976, Leibniz to the Abbé Conti, 29th March 1716: 304-314, 312. This letter continued, '..excepté dans celle q'il a écrit le 24 d'Octobre 1676, où il n'en a parlé que par enigme; et la solution de cette enigme, qu'il n'a donnée que dix ans après, dit quelque chose, mais elle ne dit pas tout ce qu'on pourroit demander.' See also Leibniz' letter to Conti of 9th April 1716 in *Corr.*, 6, p.307.
[22] The October letter cited π to 16 places, 14 of them correct.

Newton described how fractional indices and logarithmic terms could be integrated using these infinite series. In reply to these letters, in June 1677 a letter from Leibniz to Oldenburg the Royal Society's Secretary, described the reciprocal nature of the two operations, differentiation and integration, and the Leibnizian notation of dy and dx. Hall found that, 'There was, in fact, a good deal more basic calculus in this letter of June 1677 than Newton was to make public in his *Principia* lemma or anywhere else before 1704'[26]

Integration is less of a general method than differentiation, and remains to this day more of an art; one has to find the integral by trying methods until one works. The differential technique by contrast is more or less automatic, and here bears comparison to Leibniz' new calculating engine, or his search for a universal symbolic logic, which would give an automatic procedure for evaluating certain truths. What could have been acknowledged, was that certain difficult integrals had been achieved by the English maestro. It is hard to think what else could or should have been acknowledged, by way of supporting Westfall's position. To quote Hall again, 'Leibniz had it fixed in his mind, all his life, that Newton's great expertise was in the manipulation of series – nothing else.' (Hall *Philosophers* p.66)

## 1680s: the *Arithmetica*

Newton's *Arithmetica Universalis* was based upon his lecture-notes of the 1680s as Lucasian mathematics professor, though only published in the next century (1707) by William Whiston. While nowadays forgotten, it enjoyed considerably more reprints in the eighteenth century than did the *Principia*. Not only does it contain no trace of differential calculus ('fluxions'), but it warns students against trying to compound geometry and algebra:

> 'Equations are Expression of Arithmetical Computation, and properly have no Place in Geometry. The ancients did so assiduously distinguish them from one another, that they never introduced Arithmetical Terms into Geometry. And the Moderns, by confounding both, have lost the Simplicity in which all the Elegancy of Geometry consists.'[27]

We are here reminded of his later notes for a '*Geometriae*' text of the early 1690s (unpublished) where he argued for the superiority of geometrical methods for 'obtaining direct and natural ways of solution', whereas modern mathematicians were prone to reaching less simple solutions because they are more 'devoted to algebra.'[28]

---

[23] De Morgan 1852, p.88.
[24] Hall, 1980, p.95.
[25] In Henry IV Part 1, Owen Glendower says: 'I can call up spirits from the vasty deep' and Hotspur replies, 'And so you can; and so can I; and so can any man: but do they come when you do call for them?'
[26] Hall, 1980, p.71.
[27] Isaac Newton *Universal Arithmetick* 2nd Edn 1728, trans from the Latin by Mr Raphson, p.227-8.

## The *Principia*

*There are no extant autograph manuscripts of Newton's preceding the Principia in time which could conceivably buttress the conjecture that he first worked the proofs in that book by fluxions before remoulding them in traditional form. Nor in all the many thousands of such sheets relating to the composition and revision of the Principia is there any trace of a suggestion that such papers ever did exist* – Tom Whiteside[29].

De Morgan characterized the *Principia* as 'the work of an inordinate Euclidian, constantly attempting to clothe in the forms of ancient geometry methods of proceeding which would more easily have been presented by help of algebra.'[30] In the text of the *Principia*, as Whiteside well remarked, 'the geometrical limit-increment of a variable line-segment plays a fundamental role.'[31] That is a statement about geometry, not algebra. Its never-fathomed diagrams contain a springy movement within them. Their limit-increments blur the sharp outlines of Greek geometry; but, it is a proto-calculus, and of the integral kind[32]. Its 'moments' or indefinitely small increments are always connected with geometrical diagrams.

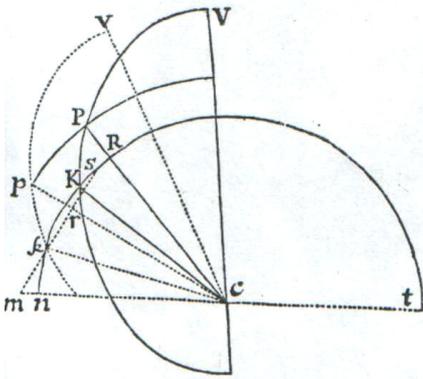

We must surely agree with De L'Hôpital that the *Principia's* author 'avoit aussi trouvé quelque chose de semblable au calcul differential…'[33] Its geometrically-structured arguments did indeed carry a 'semblance' of differential calculus. Primarily, the *Principia* used integral methods of calculus, not differential. We here recall Leibniz' comment as to how he had come to his new calculus 'by which the imagination is freed from a perpetual reference to diagrams.'

The *Principia* defined *momenta* as 'only just-born finite magnitudes' (*jamjam nascentia finitarum magnitudinem*) with the essential paradox that –

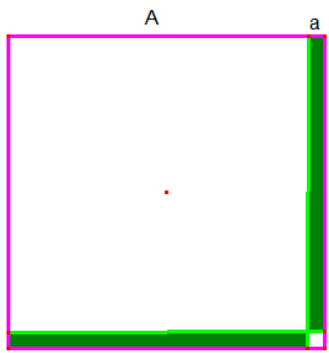

*Momenta, quam primum finitae sunt magnitudinis, desinunt esse momenta.*

– as soon as they reached a finite magnitude, they would cease to be moments. The term 'fluxion' was *not used* in the *Principia*, being first introduced in the 1693 Wallis opus, despite which the above-quoted section in its Book II has misleadingly come to be called 'the fluxions lemma.' It explained how, if a quantity A suffered a small 'mutation' 'a,' then the function $A^2$ would have the 'moment' 2Aa. If A expands a little to become A+a, then as the diagram shows its original area $A^2$ will have expanded by the shaded areas equal to 2Aa, ignoring the small corner area $a^2$. A postscript to this lemma averred that this method was that signified to Leibniz in Newton's early 1672 letter, concerning 'methodi determinandi Maximas & Minimas, ducendi Tangentes'. Neither in this Lemma nor anywhere else in the *Principia* did Newton give any such method for finding tangents or maxima and minima, yet seemingly by this comment he expected to be given credit for such (a different method for finding tangents but not maxima or minima had earlier been advocated in Newton's 1672 letter to Leibniz, viz. that of Sluse).

---

[28] D.T. Whiteside, *The Mathematical Papers of Isaac Newton,* (MP) Vol. VII, CUP, 251, *Geometriae* Preface.
[29] Tom Whiteside, *The Mathematical Principles underlying Newton's Principia* 1970, U of Glasgow, pp.9,10.
[30] De Morgan 1852, p128
[31] D.T.Whiteside, op.cit. 1960, p.10.
[32] Newton's 1687 *Principia* diagrams: p.132 Elliptical trajectory, p.109 Apse motion: Domenico Meli, *Equivalence and Priority, Newton versus Leibniz,* 1993, Oxford, p.70.
[33] De L'Hôpital, 1696, Preface: '… likewise discovered something like the Calculus differentialis.' It continues, 'But the method of Mr Leibniz's is much more easy and expeditious, on account of the Notation he uses…' (Translation from J.Fauvel and J.Gray, Ed., *The History of Mathematics – a Reader*, OUP, 1987, p.442).

According to Hall, this lemma 'explains simple differentiation (as we call it)'[34] That has to be a great exaggeration: what has been described is a recipe for ignoring power terms when dealing with small change, e.g that $(A+a)^3$ equals $A^3 + 3A^2a$ to a first approximation, the latter term being the 'moment'. No rate-of-change concept was there present nor anything resembling the concept of a differential, viz. the limiting value of a ratio between two infinitesimal quantities.

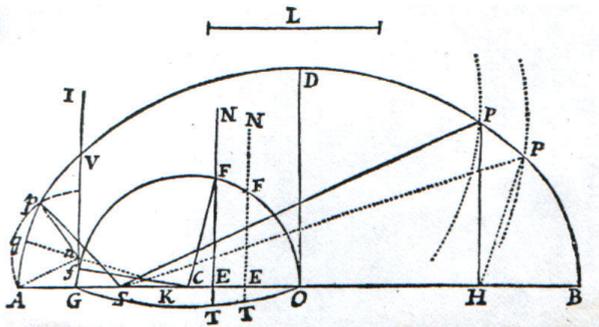

One may here accept Leibniz' view, 'It was evident from the *Principia* that Newton had not known the true Calculus of Differences in 1687, because he had avoided its use in circumstances which demanded it'[35]. To give an example, Book 3 brilliantly explained why two tides happened per day, via an inverse-*cube* law of tidal pull[36] - however it gave readers no hint whence this had been derived. With his deep intuition Newton had apprehended that the tidal pulling-force would have to vary inversely as the cube of the luminaries' distance from Earth, whereby he could explain how the Moon had a greater pulling power on the tides than the Sun, though the latter's gravity pull is always far stronger. Those concerned to popularise Newtonianism such as Halley and Gregory had to be content with merely sketching out an argument in rough qualitative terms, because they had no way of following how that inverse-cube was supposed to work. The inverse-cube law of tidal pull is derivable very readily - once differential calculus is available.

Westfall's biography routinely alludes to Newton's 'g' computation[37], an illusory concept which has featured heavily in popular accounts. The actual computation, well described by Dana Densmore as 'perhaps the most thrilling demonstration in the *Principia*,'[38] forged a link between Heaven and Earth, by comparing the 27-day lunar sidereal orbit-period with the fall of objects on Earth: the distance which the Moon would fall Earthwards in one minute 'if stopped' was compared with the distance which objects on earth fall in one second. It used *no concept of acceleration* or even velocity.[39] This was achieved in the year 1685.[40] Propositions 32-39 of the *Principia's* Book I (third edition) deal with the acceleration of falling bodies, but

---

do so in a geometrical manner, and without assigning a numerical value to acceleration. Popular writers have interpreted what Newton did in terms of differential concepts, such as acceleration as distance per second per second, but we need to apprehend that such interpretations only appeared in the mid-eighteenth century – using the Leibnizian calculus - in France and Germany.

Newton's second 'Law of Motion' affirmed that 'change of motion is proportional to the motive force impressed'[41], Its seven Latin words "*Mutationem motus proportionalem esse vi motrici impressae*" are a statement about impulse, however no-one could count the number of commentators who have read this as a statement about rate of change[42], absent though this is from the text[43]: as if it lead to or even was the differential law connecting force and acceleration invented by Euler in 1737[44]. How many generations of children have been taught this law as F=MA? One would be hard pressed to find an English textbook that gives credit to Euler for his fundamental discovery, owing to their near-unanimity in ascribing it to Newton! Calculus appears in a half-born condition in the *Principia*'s text, not yet formulated. As Bernard Cohen concluded, Newton 'did not ever write' a fluxional equation such as F = kms" or F = kmv̇ as we would.'[45]

Many have been tempted to presume that Newton was in the 1690s using his newly-defined fluxional terms ẋ, ẏ for speed, distance and acceleration computations. Whiteside has well described Newton at this time as 'caged by its [the *Principia's*] very brilliance and originality, unable to transcend its mental confines' (MP, 7: xi) which may help us to appreciate why he did *not* in fact start sending 'fluxional equations' to David Gregory or anyone else. Instead, he had a nervous breakdown, then next year focused his mind upon lunar studies, which became through the nineties his last great scientific endeavor.[46]

The Italian scholar Nicollo Giucciardini is continually describing the *Principia's* arguments as 'fluxional,' with fluxional equations etc.: 'His only coherent policy was that of not disclosing in print the analytical parts of his demonstrations.'[47] Newton composed this work in just over two years, so he worked rather quickly, and we should surely agree with Hall and Whiteside that the extensive notes which Whiteside has combed through appear not to include prior 'analytical' or 'fluxional' reasonings.

No-one was closer to Newton's mathematical work than the Scotsman David Gregory, and Guicciardini has him 'interested in analytical fluxional versions of the *Principia's* demonstrations'.[48] Was he? I checked Gregory's extensive *Notae* held in the Royal Society's library, which are a commentary on the *Principia*. They have infinite series expansions, and some 'quadratures,' using 'moments' defined by geometrical diagrams - but not, I suggest, any fluxional equations. In May 1693, did Newton send to Gregory 'the solution of the pertinent fluxional equation'?[49] I suggest not, although this may be to some extent a semantic issue. Gregory made notes after his formative visit to Newton over 4-8 May 1694, and these state:

> *The problem of quadrature and the inverse method of tangents includes the whole of more advanced geometry*[50]

---

[41] Sir Isaac Newton's *Principia*, Motte translation, Berkeley, California U.P, 1934, p.13.
[42] See Herman Erlichson, 'Evidence that Newton used the calculus to discover some propositions in his "Principia," *Centaurus* 1997, 39: 253-266, p.263, for the erroneous claim that the Second Law concerned rate of change of motion: using, moreover, dv/dt!
[43] E.J. Dijksterhuis, *The Mechanisation of the World Picture*, 1961, Princeton, Princeton U.P. p.471, para 303.
[44] C.Truesdall 'A program towards rediscovering the rational mechanics of the Age of Reason', *Archive for History of Exact Sciences*, 1960, 1, p.23.
[45] I Bernard Cohen, *Introduction to Newton's Principia* 1971 CUP, p.166.
[46] N.K., [Newton's forgotten Lunar Theory](#), Green Lion books, 2000.
[47] Guiccardini, N., *Reading the Principia* CUP 1999, p.95
[48] Ibid,, p184.

Futurity would not see it that way at all, rather the contrary. Gregory hoped that his extensive commentary on the *Principia* would be published (it never was) and he enumerated in July 1694 things that needed to be added in. These included:

> *whatever things hitherto as yet unpublished are kept close by geometers or sons of the Art, these are to be inserted in their proper places, as the Differential calculus of Leibniz and other matters.'*[51]

Gregory appears as a major conduit whereby the Leibnizian calculus was transposed into the British fluxional notation[52] and his 1695 tract *Isaaci Newtoni Methodus Fluxionum, ubi Calculus Differentialis Libnitii…*[53] compared the continental and British methods. Eagles concluded that 'Gregory did not so much absorb the new concept of 'fluxion' but redefined it in terms of the 'moments' with which he was used to dealing' – the latter being the Leibnizian differential.

## John Wallis' *Algebra*, 1685 & 1693

*'The letters to Wallis in 1692 … [were] the first significant announcement to the world at large of the power of Newton's fluxional method.'* - Tom Whiteside[54]

John Wallis, the Savilian mathematics lecturer of Oxford, was known for his assertive promoting of English mathematics against that of the Continent. The first edition of his *Algebra* of 1685 had contained no hint of any new method of calculus.[55] The British mathematical innovations which he there highlighted mainly concerned infinite series expansions. Newton's unpublished 1669 manuscript *De Analysi* was cited for a couple of features: its nomenclature for power terms, e.g. writing $a^{1/2}$ for the square root, and also a method of approximate solution to equations, which was not iterative, nor had any fluxional aspect to it[56]. Twentieth-century scholars tended to assume that *De Analysi* contained such a differential method, to the extent that the mere question of whether Leibniz may have seen it almost suffices to resolve the priority dispute, or to convict him of dishonesty for not admitting what he (supposedly) gained from it. But, *no-one* in the seventeenth century, amongst those who saw the manuscript as circulated by John Collins, is on record as attributing to it such a significance[57]. Only towards the end of the first decade of the 18$^{th}$ century, when it was becoming evident that the invention of calculus was the biggest thing since the Arabs had found zero, did such comments start to appear. Sagely cautioned Tom Whiteside,

> *'Be ever wary of accepting uncritically at face value Newton's retrospective assessments of the nature and sequence of his scientific achievement'*

---

[49] Ibid., p.183
[50] Memoranda by David Gregory, *Corr.* III, 318.
[51] *Corr*. III:386.
[52] Christina Eagles, 'The Mathematical Work of David Gregory 1659-1708', Edinburgh PhD 1977, p.405-422.
[53] ULC Add. 9597.9.3 &4 (two copies), original at St. Andrew's.
[54] 'Essay Review of Correspondence, Vol. III' History of Science 1962, 1: 97.
[55] John Wallis, A Treatise on Algebra, 1685, p.338.
[56] N.K., 'John Simpson and 'Newton's Method of Approximation': an Enduring Myth', *BJHS*, 1992, 25, 347-54. An iterative approximate-solution method given in Joseph Raphson's 1697 book (see below) became known much later on as 'Newton's method', and one finds it cited as being present in *De Analysi,* (eg, J. Pepper, Ch.3 of *Let Newton Be! A new perspective on his Life and Works,* Ed. Fauvel et. al., OUP 1988, 73; D. M. Burton, *The History of Mathematics, an Introduction*, 1986, 408). Thus a rich texture of mythology surrounds these events.

Adding his all-important conclusion:

> *'As Newton himself never forgot (though he tried cleverly to conceal it by various sly turns of phrase from his contemporaries) the standard dot-notation for fluxions was invented by him only in mid-December 1691 ...58.'*

Newton communicated his fluxional method in a letter of September 1692, now lost, to Wallis [59]. This method concerned 'flowing quantities and their fluxions,' and gave the dot-notation for time-derivatives: thus $x^3$ would yield the fluxional term $3x^2\dot{x}$.[60] Wallis published this method in 1693[61],[62], but without alluding to the letter he had just received from Newton. Instead, he alleged that certain letters Newton had sent to Oldenburg for Leibniz in 1676 contained this method - as was plainly not the case[63] - whereby seeds of far-reaching confusion were sown. This text was 'really Newton's first publication on the subject [of fluxions]' as Augustus de Morgan correctly called it,[64] or 'the first public account of the intertwined methods of fluxions and infinite series'[65], to quote Whiteside. The shy debut of fluxions on the English stage was a method of implicit differentiation for time-dependent functions – and *did not find* the gradient of curves.[66] This 1692 fluxional method reappeared unchanged in *De Quadratura*.

## Bernoulli's Challenge

Johann Bernoulli posed a challenge problem in the June 1696 issue of *Acta Eruditorum,* designed to sort out which mathematicians in Europe had mastery of the 'golden equations' as he called them, of the new calculus. For he was very confident, that no-one could solve his challenge, concerning changing rate of motion along a curve, without them:

> *Given two points A and B in a vertical plane, what is the curve traced out by a point acted on only by gravity, which starts at A and reaches B in the shortest time?*

---

[57] Newton's letter to Oldenburg of 24 October 1676 alluded to this manuscript, as having been communicated by Mr Barrow to Mr Collins, and averred that its quadratures could find 'areas and lengths of all curves, and the surfaces and volumes of solids from given right lines.' (*Corr.* I, p.133.)

[58] Tom Whiteside, *The Mathematical Principles underlying Newton's Principia* 1970, U of Glasgow, pp.9,10.

[59] *Corr.,* III, 222-8: 'Newton's Method of Fluxions', a transcription of what Newton presumably sent, dated 17$^{th}$ September, 1692.

[60] Newton's example was $x^3 – xyy + aaz – b^3 = 0$ with x, y and z as variables, from which he obtained the fluxional equation $3\dot{x} x^2 – \dot{x} yy – 2xy\dot{y} + aa\dot{z} = 0$. This same example was re-used *in De Quadratura (Optics*, 1704, p172).

[61] J.Wallis, *Opera* II, 1693, *A Treatise on Algebra*, 2$^{nd}$ Ed.,391-6; partly reproduced in *Corr.* III 1961, 220-1.

[62] Leibniz commented, 'but that he [Newton] was occupied with a calculus so similar to the differential calculus, was not known to me until the appearance of the first of a two-part work by Wallis.' (*Responsio* to Fatio de Duller)

[63] Wallis re-presented this message in a more polemical form in the Preface to Volume I of his *Opera* published two years later in 1695: that Newton's fluxional method given in volume II was both equivalent to that of Leibniz and published a decade earlier.

[64] Augustus de Morgan, 1852, p.92.

[65] MP, VII, General Intro, p.xvii.

[66] Hall (1980, p.258) claimed that Newton's fluxional letter sent to Wallis in 1692 described 'its use in problems of tangents, maxima and minima.' It really did not do that.

The six-month deadline elapsed with no English response, so Leibniz persuaded him to extend the deadline, and then Bernoulli posted a copy to Newton: after all there had been rumors that Newton was in possession of the new calculus and would this not resolve the matter? John Wallis in Oxford, David Gregory in Scotland, Pierre Varignon in Paris – they all struggled in vain, and had to admit defeat.[67] Only four answers were received: from Leibniz, De L'Hopital, Bernoulli's brother Jacob, and Newton.

It is pertinent to quote the view of Pierre Rémond de Montfort, concerning Newton's solution:

> *'On trouve en 1697 une solution de M Newton du probleme de la plus viste descente, main comme il n y a point d'analyse & qu'on ne scait point la route qu'il a suivi…'*[68]

The *purpose* of that challenge question was to show how the new calculus was required to find a curve, to ascertain which mathematicians in Europe were in possession of its 'golden theorems' (to quote Bernoulli) capable of demonstrating these things. Newton sent off on 30$^{th}$ January a two-sentence solution, which *presupposed* that the catenary curve was to be used, indeed it was little more than the joining up of two points using a catenary; then a month later a slightly fuller version was published in the *Phil Trans*.[69] The question of the path of quickest descent of a falling body could hardly have been more pertinent to Newton's work. Bernoulli's May 1697 published solution[70] *demonstrated* that a cycloid curve resulted from such a path of quickest descent, whereas Newton's had merely submitted an 'undemonstrated construction of its cycloidal curve,' to quote Whiteside.[71] It was no mystery that a cycloid was the curve traced out by the edge of a rolling wheel. There was a second problem also to be solved in this challenge and Newton's solutions to both were purely geometrical – he *did not use* his 'fluxions' to solve them. A consequence of his not having utilized the calculus, in de Montfort's view, was that one could *not follow* what Newton had done.

Whiteside conjectured that Newton had actually used his maximum/minimum infinitesimal method on this problem (MP 8:6, n.12) – maybe or maybe not, but we only see a Newtonian text using fluxions to demonstrate this matter in March 1700 – a year *after* Fatio Duiller had published a book on the subject (cited in the list of fluxion publications, below) and three years after the Leibnizian solutions were being debated in German and French journals - and even that Newtonian text remained unpublished (MP 8:86, Appendix 2). What eventually appears as a Newtonian solution demonstrating that the curve was a cycloid emerges only though discussions with David Gregory *after* Newton has sent off his response to the challenge problem (thereby rescuing England's mathematical reputation) and as a kind of echo of the solutions given in *Acta Eruditorum*. This ended up years later in Hayes' *Fluxions* of 1704 as a worked example - with no acknowledgement of its Continental origin. Nothing could be more pertinent to gravity theory than this problem of a path of descent, a challenge problem emerging from what Leibniz had named 'dynamics' that was transforming mathematics in Europe.[72] The stark fact is that Newton did not ever use his 'fluxional' methods in any public context, he preferred geometry.

---

[67] MP 8, p.5
[68] Letter, Monmort to Brook Taylor, 7.12.1718,. *Corr.* VII pp.21-22.
[69] Newton to Montague Corr., 4, p226, then PT 97,19, 384-9.
[70] Bernoulli, *Acta Eruditorum* May 1697 206-217; also a solution by De L'Hopital 217-223.
[71] MP 8:10. As R. Woodhouse observed (*Treatise on … the Calculus of Variations* Cambridge 1810), Newton had merely submitted 'without proof… a method of describing the cycloid,' p.150.

## British Fluxion Publications

Here is the sequence of British publications concerning the Newtonian method of fluxions, from 1695 up to 1704 when Newton's *De Quadratura* was published. We here include by way of contrast three (in blue) by John Craig that alluded only to the Leibnizian method:

\* 1685 John Craig *Methodus figurarum ... quadraturas determinandi* inspired by Leibniz' *Acta E.* article of 1684 (p.27), plus also Barrow and Gregory. De L'Hopital alluded to its part II. It concerned finding areas (ie integration), with little hint as yet of fluxional method.
\* 1693 J. Craig *Tractatus Mathematicu de figurarum Curviliearum Quadraturis* (alludes to Leibniz, *Acta E.*1684, also to Barrow). Basic Leibnizian methods, 76pp.
\* 1695 Abraham de Moivre, *Specimia ... Doctrinae Fluxionum* Philosophical Transactions 19, pp.52-7. Alludes to *Principia's* 'Fluxions lemma' and to Craig; uses both the dx and the x fluxional notation. Is concerned with calculating solids of rotation and centres of gravity, i.e integration.
\* 1699 Fatio Duillerii *Lineae Brevissimi Descensus* 24pp.
\* 1697 Joseph Raphson *Analysis aequationum universalis*.
\* 1702 John Harris *A New Short Treatise of Algebra* with a last chapter, 'of Fluxions' 21 pp.
\* 1702 De Moivre, *Methodus Quadrandi Genera Quaedam Curvarum* PT 1113-27 uses mainly Leibnizian differentials, but Newtonian fluxions in one section.
\* 1703 George Cheyne, *Fluxionum Methodus Inversa* (infinite series, no diagrams, no max/min found) alludes on p2 to Wallis' *Algebra*, p.293; 128 pp. Alludes to Craig, 1685. (De Moivre PT230 series expansions)
\* 1703 J. Craig Specimen Methodi Generalis determinandi Figurarum Curvarum, uses integral and differential signs: PT 23, 1346-60 'Methodum Calculi Differentialis,' no fluxions.
\* 1703 De Moivre, *Animadversions in Georgii Cheyneai Tractatum de Fluxionum* 129pp
\* 1704 Charles Hayes *A Treatise of Fluxions* 315 pp.

These works employing fluxions alluded to the two Newtonian sources viz. the 'fluxions lemma' of the *Principia*'s Book II, and the 'fluxions' section in John Wallis's 1693 *Algebra*, 2$^{nd}$ Edition, which quoted some Newtonian letters. John Craig in his 1685 opus (above) alluded to the just-published Leibnizian calculus, as likewise did his next book on the subject. In the first of these books, Craig says he has been permitted to inspect Newton's manuscripts, indicating that Newton had helped him with certain infinite series expansions. De Moivre's *Philosphical Transactions* article of 1695 alluded to the earlier Cheyne opus. Newton's *De Quadratura* alluded to the1703 George Cheyne volume, as having published 'his' method, and as the stimulus for Newton finally going into print.[73]

These five above-cited British fluxion-textbooks beginning in 1695 each reiterated the product and quotient rules for differentiation, as given in Leibniz' 1684 article in *Acta Eruditorum.* This had affirmed, using the Leibnizian differentials, that:

$$d(xy) = ydx + xdy, \text{ and that } d(x/y) = (ydx - xdy)/yy$$

---

[72] In 'Newton and Leibniz's dynamics,' P. Costabel (*The Texas Quarterly,* 10:119-26) contrasted the 'dynamics' of Leibniz with the 'rational mechanics' of the *Principia*.
[73] Whiteside found texts for *De Quadratura* composed around 1693 i.e the time of the letters to Wallis: MP VIII.

The British publications replaced dx by the Newtonian fluxion notation ẋ , and dy by ẏ . These were equivalent, alluding to an indefinitely small quantity.[74] Leibniz' article had explained how maxima and minima of curves were found by putting the differential dy equal to zero, and also how points of inflexion, where a curve turns from being convex to concave, were found by putting the second differential 'ddx' equal to zero. *Only* Charles Hayes' 1704 book in this series gave a fluxion-based explanation for such maxima and minima of curves - taking five pages for this – as likewise he did for curve points of inflexion.

John Craig's *Method of determining the Quadratures of figures* of 1685 gave Newton's circle-squaring formula, i.e the quadrature of a circle to find its area. Prior to 1696 when de L'Hôpital's textbook *L'Analyse d'Infiniments Petits* expounded the Leibniz-Bernoulli method, only one rather slight Newtonian text had been published, that by Abraham de Moivre: this was a riposte to the earlier Craig opus, using the same equation as had Craig but showing how the implicit differentiation process worked using fluxions instead of differentials. His text used *both* the dx and the ẋ fluxional notation.

Three years before Newton's *Principia,* Craig's text opened with the acknowledgment: 'that very famous person [Leibniz] shows a neat way of finding tangents even though irrational terms are as deeply involved as possible in the equations.' Newton read through the manuscript and made some suggestions, before it was printed, without objecting to a text where only Leibniz' contribution was acknowledged[75]. Craig's attempt to interest his fellow-countrymen in the Leibnizian calculus proved futile, the only sequel to his opus being his own second edition in 1693, which said: 'I freely acknowledge that the differential calculus of Leibniz has given me so much assistance in discovering these things that without it I could hardly have pursued the subject with the faculty I desired; how greatly the very celebrated discoverer of it has advanced the solid and sublime art of geometry by this one most noble discovery cannot be unknown to the most skilled geometers of this age, and this treatise now following will sufficiently indicate how remarkable its usefulness has been in discovering the quadratures of figures.'[76]

The above British publications over the nine-year period 1695-1704 utilised the new Newtonian concept of 'moments' labelling them as fluxions. Thus, in 1702, John Harris wrote, 'These infinitely small Increments or Decrements, our Incomparable Mr Newton, calls very properly by this name of fluxions,' (p.115) and likewise in 1704, Charles Hayes wrote: 'the infinitely little Increment or Decrement is call'd the fluxion' (p.1) In 1711, John Hayes' *Lexicon Technicum* had a twelve page section on fluxions which started with such a definition: 'the fluxion of a quantity is its Increase or Decrease infinitely small.' These British publications do not resolve the dilemma, as to whether a definite value could be assigned to an indefinitely little quantity, and *for none of them* does a fluxion signify a rate of change (dx/dt in differential form). They remain as little quantities, little enough that squares and higher powers can be ignored.

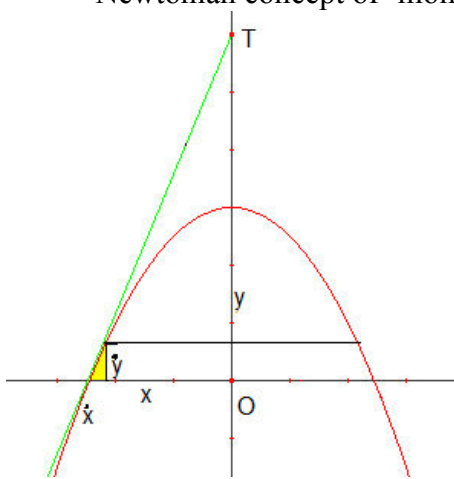

---

[74] We here exclude the Fatio Duiller book, as focused upon only one curve, the 'cycloid'.
[75] Hall 1980, p.79.
[76] John Craig, *Tractatus Mathematicus de Figurarum Curvilinearum Quadraturis* 1693, London, p.1.

Leibniz's original 1684 announcement of the method had suffered from this same dilemma, as to whether a definite value could be assigned to an indefinitely small increment. This central problem was not resolved until 1696 when de L'Hôpital started writing dy/dx: one indefinitely small increment divided by another can indeed have a definite value, that of a limit towards which it moves ever closer.

Did any of the above sequence of British publications 1695-1704 find the gradient of a curve? In 1702, Harris 'found' a tangent, to a parabola. His tangent to the parabola py = xx cuts the axes at x and y. The fluxional form of the equation is $p\dot{y}$ = $2x\dot{x}$. By similar triangles, $\dot{x}$ : $\dot{y}$ = x : OT, if the tangent cuts the y-axis at T. Substituting in gives him OT = 2y, so that the tangent intersects the y-axis at twice the value 'y' where the parabola cuts it. QED, that was it! That 'similar triangles' argument is the nearest Harris got to a statement about gradient. We are rather startled that he should so brazenly have lifted this worked example from de L'Hôpital without acknowledgement, merely putting $\dot{x}$ for dx.[77]

There was however one difference, in that de L'Hôpital had written, concerning his tangent example (repeated by Harris), PT = y dy/dx. That little triangle *defining the gradient* was then 'infiniment petits' and

> *'donnera une valeur de la soutangente PT en termes entiérement connus & délivrés des différences'*

It is not going to far to describe 'une valeur de la soutangente' as the gradient of a tangent, and so the finding of this 'freed from the differences' means that it has become independent of the indefinitely-small values of the dx and dy terms. That is the all-important step, after which the equals sign can be properly used concerning quantities that are finite i.e. of knowable magnitude. The problem of how to express the infinitely small as a limit-ratio was thereby resolved, in France.

Hayes' 1704 textbook used the same format as de L'Hôpital had done for finding the tangent to a curve, just putting $\dot{x}$ instead of dx, i.e. he gave the gradient of his parabola-tangent as $\dot{x}$ / $\dot{y}$, a ratio between two infinitesimal increments (p.18). Thus a British publication first expressed the gradient of a curve - eight years after the concept had been established in France.

Raphson's *History of Fluxions* of 1715 might have been less biased had it alluded to his own earlier work of 1697. Published as a tract in 1690, *Analysis aequationum universalis,* it concerned a method that was original and iterative for the approximate solution to equations. Much later on, this came to be (mis-) called Newton's method,[78] by which time it was used with differential notation. Raphson gave long lists of rules, for the operation which later came to be called differentiation, for converting eg xxxx into 4xxx (i.e., $x^4$ to $4x^3$). His book had pages of examples for all sorts of functions - with *no hint* that a general way for doing this might exist. Published the year after De L'Hôpital's textbook, it did not find the latter to be relevant, or suggest that sources more recent than Viète were available. It did not mention the differential calculus, which should have been, or would shortly become, the essence of his method. One is here reminded of De Morgan's view, that British writers started writing about the new calculus *after* de L'Hôpital's textbook appeared in Paris.

---

[77] De L'Hopital 1696, *L'Analyse d'Infiniments Petits,* p.11 'pour trouver les Tangentes'.
[78] N.K, Newton's method of approximation, an Enduring Myth BJHS,1992, 25, 347-54.

Leibniz' comment concerning John Wallis, that he showed 'an amusing affectation of attributing everything to his own nation' (writing to Thomas Burnett), could apply to several characters here alluded to. Thus, Raphson's *History* opened with the affirmation on its first page that in 1671 Newton had obtained 'a General Method of Quadratures' - and gave as his reference page 3 of a text by David Gregory, published in 1684. That was a hard-to-obtain tract, not a published book, and many of his readers might have found some difficulty in obtaining a copy to check this claim. On that page however we find, as we have seen earlier, Gregory reprinting Newton's infinite-series formula for the area of a circle, without giving any explanation for it. Thus Raphson's was a misleading argument, sliding between meanings old and new of the word, 'quadrature.' Our analysis has greatly failed to confirm Raphson's thesis, of any British line of descent for the development of fluxions. Raphson's *History* simply omitted allusion to the great French textbooks, de L'Hopital's *L'Analyse d'Infiniments Petits* and Charles Reynau's masterly *Analyse demontrée…en employant le calcul ordinaire de l'algèbre, le calcul diferential et le calcul intégral'* (Paris 1708), the latter being sometimes described as the first textbook of the new mathematics.

The first statement of Newton's method in his own name was in his *De Quadratura* of 1704, *twenty years* after Leibniz had so published. It described his integration or 'quadrature' method as: let $dZ^{n-1} = y$ be the equation of a curve, then $d/nZ^n = t$ will be the area under it, where 'd' is a constant term.[79] This process was alluded to as 'the inverse method of fluxions,' however his terminology did not really link it to his method of fluxions or moments – that was Newton's problem. Newton's text lacked any symbol comparable to Leibniz's ∫ as the active agent, and contained not a word about tangents or maxima/minima. It reiterated the process of implicit differentiation given in Wallace eleven years earlier, using the same worked example.[80] His small change was now a fluxion: "fluxio quantitatus x est ad fluxionum quantitatus $x^n$ ut 1 ad $nx^{n-1}$" So we now have a rudimentary account of something resembling differentiation, in which the ratio of the small quantities (fluxions) reaches a limit-value.[81]

Grattain-Guinness has compared 'Principal features of the calculi of Newton and Leibniz' in a table (Fontana *History of Mathematics* pp.244-5) and for 'Basic derivative concept' he put under 'Newton' $\dot{y}/\dot{x}$ and under 'Leibniz' $dy/dx$. These are both to be doubted. Never did Newton write $\dot{y}/\dot{x}$ for a gradient or derivative function, he never published a text with $\dot{y}/\dot{x}$; while the all-important $dy/dx$, originated by the Bernoulli brothers in the early 1690s,[82] appears clearly explained in the classic de L'Hôpital textbook, written in Paris as Marquis took math lessons from Johann Bernoulli.

In his text published in John Wallis' *Algebra* of 1693, Newton first defined what he called 'fluxions'. The 'moment,' defined in the *Principia* as a very small change, was there expressed as oẋ where o is a 'quantity infinitely little' and ẋ the fluxion (x being the 'fluent' quantity i.e the variable). Years later in his 1715 *Account of the Commercium Epistolicum* Newton reaffirmed this:

> *'in his (Newton's) calculus there is but one infinitely little Quantity represented by a Symbol, the Symbol o.'*[83]

---

[79] Newton, Optics 1704, *Tractatus de Quadratura Curvarum,* p.198.
[80] N. Guiccardini, in *The Development of Newtonian Calculus in Britain, 1700-1800* (CUP 1989, 2003) averred that in *De Quadratura*: 'Given an expression in which there occurred the fluxions $\dot{x}$ and $\dot{y}$ of unknown quantities, he [Newton] could manipulate the expression in order to obtain a power series of $\dot{y}/\dot{x}$. He was then able to operate on the series…' (p.2) This is just a pipedream.
[81] Newton's *Optics,* 1704, with *De Quadratura* appended, p.169, reproduced in MP 8:122-159,128.
[82] It's not in the 1684 & 1686 articles by Leibniz, which have dy and dx separate.

That is clear enough, and we therefore need to ask, why did *no-one use* this Newtonian method? Had he used it in his *Quadratura*? The fluxion ẋ has not really been defined. We have seen how the fluxion textbooks up to 1704 do not use this definition of fluxion as neither likewise did Raphson's *History of Fluxions* of 1715. Raphson rather followed the textbooks of Harris and Hayes in his Chapter II 'on Notation' in utilising the *Principia's* definition of moments whereby Fluents A, B, C etc., had increased slightly by their fluxions a,b,c etc., 'those first differences.' (p.5). As to why this conceptual confusion should have arisen – well discussed in De Morgan's ascerbic text[84] - we might conjecture that a mathematician is forbidden from dividing both sides of an equation by zero, and Newton's method did clearly involve dividing thus by 'o' which was something indefinitely close to zero. His indefinitely little concept 'o' just did not catch on. *It wasn't how to do it.*

Most of the new methods described in the above British textbooks concern methods of integration, such as how to find the volume of revolution of a curve or the centre of gravity of a solid, for which the Leibnizian differential 'dy' – or a British equivalent - was required and the Newtonian oẋ would not have been very helpful.

## A Claim for Priority

> *It would be an iniquitous judge "who would admit anyone as a witness into his own cause" announced the preface of its Report of 1712 which examined the question of propriety in calculus. Ostensibly the work of a committee of impartial scientists, the report was a complete vindication of Newton's claim and even accused Leibniz of plagiary. In fact the whole report, sanctimonious preface included, had been written by Newton himself.*
>
> Broad and Wade, *Betrayers of the Truth*[85]

In what Rupert Hall called a 'famous and delusive passage' in the 1712 Royal Society's judgment over the calculus dispute[86], Newton averred that 'By the help of the new analysis' he had found out most of the propositions in his *Principia*; but, in order 'that the system of the heavens might be founded upon good geometry' he had re-cast it into its final geometric form:

> *And this makes it now difficult for unskillful men to see the Analysis by which those Propositions were found out.'*[87]

The inscrutability of the *Principia* was now assisting his argument! Historians were thereby thrown clean off the track for two and a half centuries. Only Tom Whiteside's thorough perusal of the Newtonian documents refuted it. Rupert Hall then alluded to Newton's claim as, 'The Fable of Fluxions.'[88] By 1712 it must have been evident that his rather personal geometrical method of reasoning was not the way of the future, owing to the triumph of calculus methods: the temptation grew for that great act of deception.

---

[83] Account, *Phil. Trans.* January 1714/5, 29, p.205; Hall 1980, p295.
[84] De Morgan, 1852, p.89.
[85] W.Broad and N.Wade, *Betrayers of the Truth, Fraud and Deceit in Science,* 1985, p.28.
[86] Hall 1980, p.229, (p.296 for text of Newton's 'Account…'); adding 'It has deceived many', concerning the above-quoted passage.
[87] The 'Account' 1712 reproduced in Math Papers, 8:598-9; *Phil. Trans.* January 1714/5, 29, p.206.

We quoted Whiteside concerning the 'intertwined methods of fluxions and infinite series,' which the Royal Society's 1712 edict muddled up so badly: maybe, by following the false comment which the ageing John Wallis inscribed in 1695, in the Introduction to volume one of his *Opera,* concerning the early Newton letters. As Leibniz remarked in 1713 upon the Royal Society's judgment:

> *'They have said very little about the calculus; instead, every other page is made up of what they call infinite series. Such things were first obtained by Nicholas Mercator of Holstein, who first obtained them by the process of division, and Newton gave a more general form by extraction of roots. This is certainly a useful discovery … but it has nothing at all to do with the differential calculus.'*[89]

The Royal Society's 'Account' had two pages at the end in English. Of Newton's letter to John Collins of 10[th] December 1672 - mainly about telescope design, but with one page in effect giving Sluse's about-to-be-published recipe for finding tangents to curves, which Newton averred was his own – it stated brusquely '… in which letter the Method of Fluxions was sufficiently described to any intelligent person.' Its next page then concluded, 'For which reasons, we reckon Mr Newton the first inventor,' the 'we' being anonymous. A real nadir in scientific morality was reached when Edmond Halley acted as secretary of its clandestine group of supposed editors, for the three months of its existence, when it pronounced its judgement against the longest-standing (39 years) and most distinguished foreign member of the Society without even consulting him. This was the founder of the *Acta Eruditorum,* also of the Berlin Academy of Sciences, as well as being its President. He was also a Baron because of his archives research at Hannover: eminent in history, law, philosophy, religion, mathematics and theoretical mechanics, he shone as Europe's foremost intellectual – and he had appealed to the Royal society for Justice! Such ethics by the Savilian Professor of Mathematics must surely invite comparison with Halley's publication of his 'pirate' edition of Flamsteed's *Historia Coelestis* using the Astronomer-Royal's stolen life-work, in that same year.[90]

## Conclusion

The question is not whether you or I, or someone else, reckons they can find some germ of the differential calculus in Newton's early works: but, whether the historical characters did so. From a fairly thorough survey we can answer that they seem to have done so only right at the end of the 17[th] century, and for motives of national pride or maybe adulation of Newton. Whereas, immediately following Leibniz' two publications in 1684 and 1686, *two* British math texts were published, using and extolling this new method (by John Craig) and then the Bernoulli brothers in Switzerland managed to decrypt Leibniz' articles (they complained about their obscurity and brevity) and that really set the ball rolling. In vain we have sought for a British source for what futurity came to call, the differential calculus. British textbooks using the Newtonian notation began in 1695, as a text on quadrature (integration), while British textbooks describing the differential method (gradients of curves, maxima and minima) began in 1702.

---

[88] Hall, 1980, 1992, p.212.
[89] Leibniz, *Historia et Origo Calculi Differentialis* 1713, published 1846, trans in Childs 1920.
[90] N.K. 'Flamsteed, John' in *Biographical encyclopaedia of Astronomers,* 2007, Ed. Hockey.

The great discovery appears as having only one source, as a mighty tree grew from a single seed, with British arguments some years later appearing as retrospective constructions.[91] For example in 1700 Pierre Varignon wrote three epoch-making equations, v = ds/dt, y (force) = dv/dt, and then later that year y = - vdv/dr for centrifugal forces, with encouragement from Leibniz[92]. We nowadays can hardly think of these concepts without these Varignon equations, and yet can one find a single British equivalent to them published, or indeed any British allusion to them, in the following decade?[93] The non-copying of these equations might be a consequence of meanings that were being assigned to the British fluxional notation.

Pierre Rémond de Montfort took issue with the Royal Society member Brook Taylor over this great issue. The latter had avowed, 'I always took Sir Isaac Newton, not only for the inventor, but also for the greatest master of it.' Montfort writing in 1718 insisted that this was 'an error of fact:'

> *"I shall not examine here the rights and wrongs of Messrs Newton and Leibniz to the first inventors of the differential and integral calculus. I shall report to you when you would like the details of the reflexions that a long and serious examination have furnished me, and I hope you will not be unhappy with them. I want only to have you note that it is untenable to say that Leibniz and the brothers Bernoulli are not the true and almost unique promoters of these calculi. Here is my opinion; you may be the judge of it. It is they and they alone who have taught us the rules of differentiating and integrating, the manner of finding by these calculi the tangents of curves, their turning points and points of inflection, their greatest and their least ordinates, evolutes, caustics of reflection and refraction, quadratures of curves, centres of gravity, centres of oscillation, and of percussion, problems of inverse method of tangents such as this for example (which lent so much admiration to Mr Huygens in 1693): to find the curve of which the tangent is to the part intercepted on the axis in a given ratio. It was they who first expressed mechanical curves by means of equations, who taught us to separate the variables in differential equations, to reduce their dimensions, and to construct them by means of logarithms, or by means of rectification of curves when that is possible; and who by pretty and numerous applications of these calculi to the most difficult problems of mechanics, such as those of the catenary, the sail, the elastic, the curve of quickest descent, the paracentric, have put us and our successors on the path to the most profound discoveries. Those are facts not to be contradicted. To convince oneself of them it suffices to open the journals of Leibniz [i.e. the Acta Eruditorum.]"*[94]

---

[91] Newton wrote to Gregory 14 July 1694 explaining a Book II *Principia* lemma concerning the shape of a boat-prow offering minimal resistance to water, by equating a fluxional term to zero (*Corr*. 3:380-2) and this solution became a showpiece of the new calculus, eg Hayes 1704, 147-152.
[92] E.J.Aiton, 'The Vortex Theory of Planetary Motions', *Annals of Science,* 1958, 14, 157-172, 161.
[93] See, eg, John Keill's paper on central forces, 'Epistola … de Legibus Virium Centripetarum', PT 1708, 26:174-8. This text is Newtonian in that its 'fluxions' are small increments, always defined by geometrical diagrams. Giucciardini, 1999, p.227. A comparison with the major French textbook, Charles Reynau's *'Analyse demontrée…en employant le calcul ordinaire de l'algèbre, le calcul diferential et le calcul intégral'* (Paris 1708), no copies of which exist in the UK, would be of interest.
[94] Letter, Monmort to Brook Taylor, 7.12.1718,. *Corr.* VII pp.21-22.

We have to conclude firstly, that Newton *did not ever discernably use* 'his' method of fluxions in a public context. Questions of whether he had somehow implicitly used it, prior to the inscrutable geometric format in which his proofs were cast, are always going to remain indeterminate. Second, he defined publicly what he meant by 'fluxions' in 1693 and 1704, whereas British textbook writers and mathematicians seem not to have used this, preferring instead the Leibnizian differential but using the $\dot{x}$ notation. The blossoming of the new calculus *came from* a Continental source and that which appeared in England was an echo, always lagging by some years behind. Future development of British mathematics was impeded by the use of Newtonian fluxions for just over one century before it was discarded, having not really contributed anything but being testimony to the semi-divine status accorded to Newton.

In science, as Robert K. Merton aptly put it, nothing belongs to you until you give it away.[95] The production of private manuscripts from forty years earlier, where the author has neglected to use the methods there given because he preferred a different way, lacks any genuine relevance as regards priority of invention. Through an organic sequence of publications starting in 1671, Leibniz, the Bernoulli brothers and de L'Hôpital ushered the new method into the world. Newton's letters to Leibniz showed how certain difficult integrations were to be accomplished and Leibniz learned from them, but we should demur at seeing these as evidence for co-invention.

---

[95] Michael Mahoney, alluding to de Montfort's letter: 'Book review of D.T.Whiteside, Math Papers Vol. VIII', *Isis* 1984, 75, 366-372.